\documentclass[12pt,a4paper]{amsart}
\usepackage{euscript,amsfonts,amssymb,amsmath,amscd}

\usepackage[T2A]{fontenc}

\usepackage[utf8]{inputenc}        

\usepackage[english,russian]{babel}
\addto\captionsrussian{}
\addto\captionsrussian{}

\newcommand{\op}{\operatorname}
\newcommand{\m}{\mathbb}

\newcommand{\Z}{\mathcal Z}
\newcommand{\Arrow}{\longrightarrow}

\theoremstyle{plain}
\newtheorem{theorem}{Theorem}[section]

\newtheorem{lemma}[theorem]{Lemma}
\newtheorem{proposition}[theorem]{Proposition}
\newtheorem{corollary}[theorem]{Corollary}

\theoremstyle{definition}
\newtheorem{definition}[theorem]{Definition}
\newtheorem{remark}[theorem]{Remark}
\newtheorem{example}[theorem]{Example}

\begin{document}

\title{
Embedding theorems for quasitoric manifolds
}




\author{Victor Buchstaber, Andrey Kustarev}
\subjclass[2010]{Primary 57S15; Secondary 57R20}
\keywords{Quasitoric manifolds, torus actions, equivariant embeddings, simple polytopes}

\begin{abstract}
The study of embeddings of smooth manifolds into Euclidean and projective spaces has been for a long time an important area in topology. In this paper we obtain improvements of classical results on embeddings of smooth manifolds, focusing on the~case of quasitoric manifolds. We give explicit constructions of equivariant embeddings of a quasitoric manifold described by combinatorial data $(P,\Lambda)$ into Euclidean and complex projective space. This construction provides effective bounds on the dimension of the equivariant embedding.
\end{abstract}

\maketitle

\section{Preface}


The starting point of our research are the following classical results.\\

{\bf Theorem (Mostow-Palais).} Let $M$ be a compact smooth manifold with a smooth action of compact Lie group $G$. Then there exists a smooth embedding $M\to\m R^N$ equivariant with respect to a~linear representation $G\to GL(N, \m R)$.\\


{\bf Theorem (Kodaira)} Let $M$ be a compact complex manifold with a positive holomorphic linear bundle (for example, $M$ possesses a~rational K\"ahler form). Then there exists a complex-analytic embedding $M\to\m CP^N$.\\

{\bf Theorem (Gromov-Tishler)}. Let $M$ be a compact symplectic manifold with an integral symplectic form $\omega$. Then there exists a~symplectic embedding of $M$ to $\m CP^N$ with the standard symplectic form.\\

The central subject of this paper are theorems on equivariant embeddings of quasitoric manifolds defined by combinatorial data. One of our tasks is to improve the classical theorems in the case when combinatorial data defines the corresponding structure on the underlying quasitoric manifold.



The key object of study of toric topology \cite{newbook} is a moment-angle manifold $\Z_P = \Z_P^{m+n}$ endowed with canonical projection map $\rho_P\colon\Z_P\to P$, where $P = P^n\subset \m R^n$ is a simple $n$-dimensional polytope with $m$ facets. The manifold $\Z_P$ may be realised as the smooth submanifold in $\m C^m$ invariant under the standard action of torus $\m T^m$ on $\m C^m$. The equivariant embedding $\Z_P\hookrightarrow \m C^m$ induces the map of the orbit spaces  $i_P\colon P\to\m R^m_{\geqslant 0}$. The manifold $\Z_P$ is a complete intersection of $(m-n)$ real quadratic hypersurfaces in $\m C^m$. For polytopes $P$ admitting {\it characteristic function} $l\colon \m T^m\to\m T^n$ there exists an $(m-n)$-dimensional toric subgroup $K\subset \m T^m$ acting freely on $\Z_P$. The orbit space of this action $\Z_P/K = M = M^{2n}$ is defined to be the {\it quasitoric manifold} $M(P, \Lambda)$. Here the polytope $P$ is a~subset in $\m R^n$ given by the system of inequalities $A_P x+b_P\geqslant 0$ and $\Lambda$ is an integral $(n\times m)$-matrix satisfying the {\it independence condition}. When the polytope $P$ is fixed, there is a one-to-one correspondence between $(m-n)$-dimensional toric subgroups $K$ acting freely on $\Z_P$ and matrices $\Lambda$ satisfying independence condition and defined up to the multiplication by an~element of $GL(n,\m Z)$.   The canonical projection map $\rho_P\colon\Z_P\to P$ splits into the composition $\Z_P\to M\to P$ and the projection $\pi\colon M\to P$ is called a moment map.

Consider spaces $\m R^n\times\m C^{q}$ and $\m R^n\times\m CP^{q - 1}$ with the trivial action of $\m T^q$ on $\m R^n$ and the standard action on $\m C^q$ and $\m CP^{q - 1}$. We will provide an explicit construction of an equivariant embedding ${M(P,\Lambda)\to\m R^n\times\m C^{q_1}}$ for some $q_1$ (see theorem \ref{main}). The composition map ${M \to \m R^n\times\m C^{q_1} \to\m R^n}$ coincides with moment map ${\pi\colon M\to P\subset \m R^n}$. It follows that the moment map can be extended to an equivariant embedding.

The next main result describes an equivariant embedding ${M(P,\Lambda)\to\m R^n\times\m CP^{q_2 - 1}}$ for some $q_2$ that is also an extension of the moment map $\pi\colon M\to P\subset \m R^n$. This embedding is determined by the combinatorial data $(P,\Lambda)$ and a representation $\tilde k\colon K\to\m T^1$ (see theorem \ref{projmain}). We describe the induced homomorphism of cohomology groups $H^2(\dot,\m Z)$ in terms of the representation $\tilde k\colon K\to\m T^1$. We provide an algorithm of constructing the embedding from theorems \ref{main} and \ref{projmain} starting from combinatorial data $(P, \Lambda)$ and a character $\tilde k$. Moreover, theorems \ref{main} and \ref{projmain} provide effective explicit bounds on dimensions of the equivariant embedding.

The class of non-singular toric varieties is well-known in algebraic geometry. Every toric variety coming from simple polytope is canonically associated with combinatorial data $(P,\Lambda)$ such that ${b_P\in \m Z^m}$ and ${A_P = \Lambda^T}$. We describe implications of theorems~\ref{projmain}~and~\ref{projtoric} that allow to construct the corresponding equivariant embedding $M(P, \Lambda)\to\m CP^{q_2-1}$ and point out the connection with projective embeddings of toric varieties in algebraic geometry.

The authors are grateful to T.E.Panov and M.S.Verbitsky for valuable discussions.


\section{Notation}

$\m C^m$ is the standard complex linear $m$-dimensional space endowed with the canonical basis ${e_1=(1,0,\ldots,0),\ldots},{e_m=(0,\ldots,0,1)}$;

$\m R^m\subset \m C^m$ is the standard linear space generated by $e_1, \ldots, e_m$;

$\m R^m_{\geqslant 0}\subset \m R^m$ is the positive cone i.e. the area formed by all points in $\m R^n$ with nonnegative coordinates;

$\m Z^m\subset \m R^m$ is the standard lattice generated by $e_1,\ldots,e_m$;

$\m T^m\subset \m C^m$ is the standard compact torus ${\{(t_1,\ldots,t_m)\in \m C^m\colon |t_k|=1,\, k\in [1,m]\}}$; the map $exp\colon\m R^m\to \m T^m$ given by the formula $exp(x_1,\ldots,x_m) = (e^{2\pi i x_1}, \ldots, e^{2\pi i x_m})$ induces the canonical isomorphism of $\m T^m$ and $\m R^m/\m Z^m$. The standard $k$-dimensional compact torus is denoted by $\m T^k$ and an abstract $k$-dimensional toric subgroup in the standard torus is denoted by $T^k$;

$T_I\subset \m T^m$ is a toric subgroup corresponding to an index set $I\subset [1,m]$; the set $I=\{i\}$ defines coordinate torus $T_i\subset \m T^m$;

$\rho\colon\m C^m\to\m R^m$  is the standard moment map given by the formula 
\begin{equation}
\label{eq-rho}
\rho(z_1,\ldots,z_m) = (|z_1|^2, \ldots, |z_m|^2);
\end{equation}

$s\colon\m R^m_{\geqslant 0}\to \m C^m$ is the map given by the formula 
\begin{equation}
\label{eq-s}
s(x_1,\ldots,x_m) = (\sqrt{x_1}, \ldots, \sqrt{x_m}),
\end{equation}
note that $\rho \circ s = id$.

\section{Combinatorial data and moment-angle manifolds}

A polytope $P = P^n$ of dimension $n$ with $m$ facets is defined as the set of points in $\m R^n$ satisfying $m$ inequalities:
\begin{equation}
\label{eq1}
P = \{ x\in\m R^n\colon \left<a_i, x\right> + b_i \geqslant 0\},\, i = 1 \ldots m.
\end{equation}

We assume that inequalities are not redundant. 

The number of $i$-dimensional faces of $P$ is denoted by $f_i(P)$, so ${m = f_{n-1}(P)}$. We may rewrite (\ref{eq1}) using the matrix form:
\begin{equation}
\label{eq2}
A_P x + b_P \geqslant 0,
\end{equation}
where $A_P$ is an $(m\times n)$-matrix and $b_P\in \m R^m$.

The matrix $A_P$ and the vector $b_P$ define an affine map
\begin{equation}
\label{eq2.5}
i_P\colon \m R^n \to \m R^m, \quad
i_P(x) = A_P x + b_P,
\end{equation}
and by (\ref{eq2}) we have
$i_P(P) = i_P(\m R^n) \cap \m R^m_{\geqslant 0}$. 

The set of coordinate subspaces of $\m R^m_{\geqslant 0}$ has a canonical combinatorial structure of partially ordered set. Using this structure, we may introduce canonical moment-angle coordinates in $\m C^m$. The section $s\colon\m R^m_{\geqslant 0}\to\m C^m$ defines a map $\m T^m\times \m R^m_{\geqslant 0}\to \m C^m$ by the formula $(t, x)\to t\cdot s(x)$.
This map, in turn, allows us to identify $\m C^m$ with the quotient space $\bigl(\m T^m\times \m R^m_{\geqslant 0}\bigr)/\sim$, where $(t_1,x_1)\sim (t_2,x_2)$ if and only if $x_1=x_2$ and the element $t_1t_2^{-1}$ belongs to an isotropy subgroup of the point $s(x_1) = s(x_2)$.

Since the matrix
$A_P$ has the maximum possible rank of $n$, there exists an epimorphism $\m R^m\to \m R^{m-n}$ with its kernel equal to $\op{Im}(A_P)$. Fixing a~frame in $\m R^{m-n}$ allows to define this epimorphism by the $((m - n) \times m)$-matrix $C_P = (c_{j, k})$ such that $C_P A_P = 0$ and $C_P$ has the rank $(m-n)$.

\begin{theorem}
\label{zp}
\cite{BPR07}
The set $\Z_P = \rho^{-1}(i_P(P))$, where $P$ and $\rho$ are given by (\ref{eq1}) and (\ref{eq2.5}), is a real-algebraic $(m+n)$-dimensional submanifold in $\m C^m$. The manifold $\Z_P$ is endowed with a smooth action of $\m T^m$ induced by the standard action of $\m T^m$ in $\m C^m$. Moreover, $\Z_P$ is a complete intersection of $(m-n)$ real quadrics in $\m C^m$, which are given by the formulas
\begin{equation}
\label{eq3}
\sum\limits_{k=1}^m c_{j, k} (|z_k|^2 - b_k) = 0, \, j = 1 \ldots (m-n).
\end{equation}
\end{theorem}

\begin{definition}
The manifold $\Z_P$ is called a~{\it moment-angle manifold} of the simple polytope~$P$.
\end{definition}

If $P$ is a standard $n$-simplex $\Delta^n\subset \m R^n$, then $i_P(P)$ is a~{\it regular simplex} in $\m R^{n+1}$ -- a convex hull of basis vectors $e_1,\ldots,e_{n+1}$. The manifold $\Z_P$ coincides with the standard sphere ${S^{2n+1}\subset\m C^{n+1}}$ endowed with the canonical action of the torus $\m T^{n+1}$. Therefore, moment-angle manifold $\Z_P$ may be considered as the natural generalization of the odd-dimensional sphere $S^{2n+1}\subset \m C^{n+1}$.

Denote by $F_j,\, j \in [1, m],$ the facets (that is, faces of codimension one) of the polytope $P$.  Every index set $I\subset [1, m]$ determines a toric subgroup
\begin{equation}
\label{eq4}
T_I = \prod\limits_{i\in I} T_i.
\end{equation}
We will say that an index set $I\subset [1, m]$ is {\it admissible} if the intersection
\begin{equation}
\label{eq5}
F_I = \bigcap\limits_{i\in I} F_i
\end{equation}
is not empty. In this case every face of the polytope $P$ has the form $F_I = \bigcap\limits_{i\in I} F_i$ for some admissible set $I$.

Faces of the polytope $P$ form a partialy ordered set $S(P)$ that is called a~{\it face lattice}. We define an {\it open face} $\stackrel{\circ}{F}$ as the set of points in $F\subset P$ not lying in any of smaller faces of $P$; hence, the polytope $P$ splits into a disjoint union of open faces and every point $p\in P$ lies in a unique open face of $P$. This open face will denoted by $\stackrel{\circ}{F(p)}$ and the corresponding ordinary face will be denoted by $F(p)$. Therefore, the polytope $P$ possesses a canonical structure of combinatorial cellular complex with open faces as cells; all gluing maps in this complex are injective.

Denote by $S(\m T^m)$ the set of all toric subgroups of $\m T^m$.  Then there is a function ${\chi\colon S(P)\to S(\m T^m)}$ assigning to every face $F = F_I$ the coordinate torus $T_I\subset \m T^m$. Given a point $p\in P$, we denote by $T(p)$ the toric subgroup $\chi(F(p))$.

The canonical projection map $\Z_P\to P$ will be denoted by $\rho_P$. In this case the subgroup $\chi(F)\subset \m T^m$ is the isotropy subgroup of the set $\rho_P^{-1}(F)\subset \Z_P$ and also the isotropy subgroup of any point in ${\rho_P^{-1}(\stackrel{\circ}{F})\subset \Z_P}$. 
The combinatorial structure of $\Z_P$ considered above is induced by the canonical map $\Z_P\hookrightarrow \m C^m$. In particular, the following holds.


\begin{lemma}
\label{zpsect}
There exists a section $s_{\Z_P}\colon P\to \Z_P$ of the map $\rho_P\colon\Z_P\to P$.
\end{lemma}
{\it Proof.}
The section $s_{\Z_P}\colon P\to \Z_P$ is defined by ${s\colon\m R^m_{\geqslant 0}\to \m C^m}$.
$\Box$

\section{Characteristic functions and quasitoric manifolds}
\label{quasito}

Davis and Janusckiewicz introduced toric varieties' topological counterparts in their seminal paper \cite{DJ91}. This new class of topological manifolds has been described in terms of properties shared by $\mbox{non-singular}$ projective toric varieties. To distinguish this new class of manifolds from toric varieties, the new term {\it <<quasitoric manifolds>>} has been adopted (\cite{uspekhi2000}). 

Fixing an arbitrary set $Q$ and a group $G$, there is a fundamental correspondence between the sets of following objects (\cite{BR}):
\begin{itemize}
{\item a set $X$ with a left action of $G$ and an orbit space $X/G=Q$ endowed with a section $s\colon Q\to X$.}
{\item a map $Q\to S(G)$, where $S(G)$ is a set of subgroups of $G$.}
\end{itemize}
Given an action of $G$ on $X$ such that $X/G=Q$ and a section ${s\colon Q\to X}$, one can define a function $\chi\colon Q\to S(G)$ mapping a point $q\in Q$ to a stationary subgroup $G_{s(q)}\in S(G)$ of $s(q)$. Conversely, a~map ${\chi\colon Q\to S(G)}$ defines an action of the group $G$ on a quotient set $X = ((Q\times G)/\sim)$, where $(q_1, g_1)\sim (q_2, g_2) \Leftrightarrow q_1 = q_2$ and $g_1 g_2^{-1}\in \chi(q_1)$. We obtain space $X$ with the left action of the group $G$ given by a formula $g\cdot(q, h) = (q, gh)$ and a section $s(q) = (q, e)$. The orbit space of this action is obviously identified with the set $Q$.

Following the definition in \cite{DJ91}, let us consider a topological space
$$
\tilde \Z_P = (P\times \m T^m) / \sim,
$$ 
where
$$
(p_1, t_1) \sim (p_2, t_2) \quad \Longleftrightarrow \quad p_1 = p_2 \quad\mbox{and}\quad t_1 t_2^{-1} \in T(p_1).
$$

Topological invariants of space $\tilde \Z_P$ are determined by combinatorial invariants of the underlying polytope $P$. 
As noted in (\cite{BR}), the existence of a section map is crucial in establishing the equivalence of different constructions of moment-angle manifolds and quasitoric manifolds. In the case of moment-angle manifolds $\Z_P$ the section map exists, as shown in lemma \ref{zpsect}.

\begin{lemma}
\label{zphomeo}
There is a canonical $\m T^m$-equivariant homeomorphism $\tilde \Z_P \to \Z_P$.
\end{lemma}
{\it Proof.}
We will use the section $s_{\Z_P}\colon P\to \Z_P$ constructed in lemma \ref{zpsect}. Consider a natural $\m T^m$-equivariant map $z_P\colon(P\times \m T^m)\to \Z_P$ given by the formula
$$
z_P(p, t) = (t\cdot s_{\Z_P}(p)).
$$ 
The map $z_P$ induces a homeomorphism $\tilde \Z_P\to \Z_P$.
$\Box$

\begin{definition}
The topological manifold $\tilde \Z_P$ is called a~{\it combinatorial model} of the moment-angle manifold $\Z_P$ corresponding to the polytope $P$.
\end{definition}

\begin{corollary}
For a given combinatorial type of polytope $P$ there exists a unique (up to $\m T^m$-equivariant homeomorphism) combinatorial model $\tilde \Z_P$ of the moment-angle manifold $\Z_P$.
\end{corollary}

Recall that the geometric realization of the polytope $P$ by the pair $(A_P, b_P)$ determines the {\it smooth} manifold $\Z_P$ endowed with the canonical equivariant homeomorphism $\tilde \Z_P\to \Z_P$. There is a~correspondence:
\begin{itemize}
{\item convex geometry of $P$ $\Longleftrightarrow$ geometry and differential topology of $\Z_P$}
{\item combinatorics of $P$ $\Longleftrightarrow$ algebraic topology of the combinatorial model $\tilde\Z_P$.}
\end{itemize}

It is not known whether there exist non-diffeomorphic moment-angle manifolds $\Z_{P_1}$ and $\Z_{P_2}$ with combinatorially equivalent polytopes $P_1$ and $P_2$. As shown in \cite{davispoly}, any two combinatorially equivalent simple polytopes are diffeomorphic in the category of manifolds with angles.

The study of the combinatorial model $\tilde\Z_P$ is strongly connected with an invariant $s(P)$ known as {\it Buchstaber number} (\cite{uspekhi2000}, \cite{newbook}). The number $s(P)$ is defined as the maximum possible dimension of a toric subgroup in $\m T^m$ acting freely on $\Z_P$. Since $\Z_P$ is equivariantly homeomorphic to the combinatorial model $\tilde \Z_P$, the number $s(P)$ is a combinatorial invariant of the polytope $P$. A survey of recent results about $s(P)$ may be found in \cite{Erokhovets}.

For a given simple polytope $P=P^n$ with $m$ facets the maximum possible value of $s(P$) is equal to $m-n$. In this case there exists a~subgroup $K\subset \m T^m$ isomorphic to $\m T^{m-n}$ and acting freely on $\Z_P$. The subgroup $K$ may be defined via short exact sequence of the form
\begin{equation}
\label{eq6}
1 \longrightarrow K \longrightarrow \m T^m \longrightarrow \m T^n \longrightarrow 1.
\end{equation}

Let us fix an $(m-n)$-dimensional subgroup $K\subset \m T^m$ acting freely on $\Z_P$. Fixing basis in $K \simeq \m T^{m-n}$, we obtain an integral $(m\times(m-n))$-matrix $C$ determining a monomorphism $\m T^{m-n}\to\m T^m$. The matrix $C$ is determined up to a change of basis in the torus $\m T^{m-n}$; if $D$ is an integral matrix defining basis change in $\m T^{m-n}$, the matrix $C$ is replaced by $CD$. Once the subgroup $K$ is fixed, defining an epimorphism $\m T^m\to \m T^n$ in the short exact sequence (\ref{eq6}) is equivalent to choosing the basis in the quotient group $\m T^m/K$. Let us fix some basis in the group $\m T^m/K\simeq \m T^n$. Denote by $l\colon \m T^m\to \m T^n$ the corresponding epimorphism in the short exact sequence (\ref{eq6}). Recall that an index set $I\subset [1,m]$ is admissible if the intersection of corresponding facets is not empty in the polytope $P$. Then the homomorphism $l$ satisfies the following {\it independence condition}:
\begin{itemize}
{\item $(I \mbox{ -- admissible set}) \, \Longrightarrow\, (T_I \cap \ker l = 1)$.}
\end{itemize}

The homomorphism $l$ is called a~{\it characteristic function} for the polytope $P$.
Once the bases in $\m T^m$ and $\m T^n$ are fixed, we see that the homomorphism
$l\colon \m T^m\to \m T^n$ is determined by an integral $(n\times m)$-matrix $\Lambda$. By construction, we have $\Lambda C = 0$. Moreover, independence condition for $l$ implies the following property of $\Lambda$ (that will also be called an independence condition): 

\begin{itemize}
{\item if $\Lambda_v$ is a matrix formed by columns corresponding to facets containing some vertex $v$, then $\det\Lambda_v = \pm1$.}
\end{itemize}

We obtain a short exact sequence of coordinate tori of the form:
\begin{equation}
1 \to \m T^{m-n} \stackrel{i}{\Arrow} \m T^m \stackrel{l}{\Arrow} \m T^n \to 1,
\end{equation}
where homomorphism $l$ is given by the matrix $\Lambda$ and homomorphism $i$ is an embedding $K\hookrightarrow \m T^m$ defined by the matrix $C$.

As shown in \cite{BPR07}, this construction is invertible: any simple polytope $P\subset \m R^n$ with $m$ facets and an $(n\times m)$-integral matrix $\Lambda$ satisfying independence condition determine a characteristic homomorphism $l\colon \m T^m\to \m T^n$. The kernel $K=\ker l$ acts freely on the corresponding moment-angle manifold $\Z_P$.

\begin{definition}
(\cite{BPR07}).
The quotient space $M = \Z_P/K$, endowed with a~canonical smooth structure, is called a~{\it quasitoric manifold determined by combinatorial data $(P, \Lambda)$.} 
\end{definition}


Recall that the construction of a moment-angle manifold uses only the simple polytope $P\subset \m R^n$. To define the quasitoric manifold, one needs an additional piece of information, namely, the matrix $\Lambda$ determining a characteristic function $l$.
The problem of describing simple polytopes $P$ admitting at least one characteristic function is still open.
Simple polytopes that do not admit a characteristic function are known to exist in every dimension $n>3$. Any three-dimensional simple polytope admits a characteristic function by the four-color theorem. This characteristic function is constructed via the edge graph of the dual polytope; the algorithm is quadratic by the number of vertices of the polytope (see e.g. \cite{fourcolor}).

The manifold $M$ is equipped with the induced action of compact torus $\m T^n \simeq \m T^m / K$. The action of $\m T^n$ on $M$ is {\it locally standard}: any point $x\in M$ has an open neighborhood equivariantly diffeomorphic to a region in $\m C^n$ with the standard action of $\m T^n$. This implies that the action of $\m T^n$ has only isolated fixed points on $M$. 

If we denote by $\pi_l\colon \Z_P\to M$ the canonical projection to the quotient space, then the projection map $\rho_P\colon\Z_P\to P$ splits into composition ${\rho_P = \pi\circ \pi_l}$, where $\pi\colon M\to P$ is called a moment map. This name comes from the theory of Hamiltonian torus actions on symplectic manifolds (\cite{atiyah}, \cite{delzant}, \cite{ggk}). The moment map defines a bijection between the set of fixed points of torus action on $M$ and the set of vertices of the polytope $P$.

Let $v$  = $\bigcap\limits_{i\in I} F_i$ be a vertex of the polytope $P$. The independence condition implies that the square matrix $\Lambda_{v}$ induces an automorphism of torus $\m T^n$. It follows that the matrix $\Lambda_{v}^{-1}\Lambda$ is also the matrix of a~characteristic homomorphism $l_v\colon \m T^m\to \m T^n$, whose restriction to the coordinate subtorus $T_I$ is given by the unit $(n\times n)$-matrix. Suppose that vertices of the polytope $P$ are ordered in such a way that $v=\bigcap\limits_{i\in I} F_i$ is an~intersection of first $n$ facets. Then the matrix $\Lambda_{v}^{-1}\Lambda$ of the characteristic homomorphism $l_v$ has the following form (it is called {\it reduced form of the characteristic matrix}):
\begin{equation}
\Lambda_v^{-1}\Lambda = 
\begin{pmatrix}
1 & \cdots & 0 & \lambda_{1,n+1} & \cdots & \lambda_{1,m}\\
\vdots & \ddots & \vdots & \vdots & \ddots & \vdots\\
0 & \cdots & 1 & \lambda_{n,n+1} & \cdots & \lambda_{n,m}\\
\end{pmatrix}.
\end{equation}

In this case the $(m\times (m-n))$-matrix $C$ may be written in the following form:
\begin{equation}
C=
\begin{pmatrix}
-\lambda_{1,n+1} & \cdots & -\lambda_{1,m} \\
\vdots & \ddots & \vdots \\
-\lambda_{n,n+1} & \cdots & -\lambda_{n,m} \\
1 & \cdots & 0\\
\vdots & \ddots & \vdots \\
0 & \cdots & 1\\
\end{pmatrix}.
\end{equation}

Topological space $\tilde \Z_P$ played an intermediate role in the paper \cite{DJ91} -- it was employed for constructing topological models of quasitoric manifolds $\tilde M = (P\times \m T^n) / \sim$, where:
$$
(p_1, t_1) \sim (p_2, t_2)  \quad \Longleftrightarrow \quad p_1 = p_2 \quad\mbox{and}\quad t_1 t_2^{-1} \in l(T(p_1)).
$$

As mentioned in \cite{BR}, one has to specify a section ${P\to M}$ before proving that space $\tilde M$ is homeomorphic to an abstract quasitoric manifold defined in \cite{DJ91} by axioms. In this paper we define quasitoric manifolds as the quotients of moment-angle manifolds, so the construction of a section ${P\to M}$ is straightforward.

\begin{lemma}
\label{quasisect}
The map $\pi\colon M\to P$ has a canonical continious section $s_M\colon P\to M$.
\end{lemma}
{\it Proof.}
Let $s_M = \pi_l\circ s_{\Z_P}$.
$\Box$

The polytope $P$ has a canonical structure of smooth manifold with angles; the projection map $\pi$ and section map $s$ are smooth with respect to this structure.

\begin{lemma}
\label{quasihomeo}
There exists a canonical $\m T^n$-equivariant homeomorphism $\tilde M \to M$.
\end{lemma}
{\it Proof.}
Consider a natural $\m T^n$-equivariant map $z_M\colon(P\times \m T^n)\to M$ given by the formula
$$
z_M(p, t) = (t\cdot s_M(p)).
$$ 
The map $s_M$ induces the homeomorphism $\tilde M\to M$.
$\Box$

We see that the topological manifold $\tilde \Z_P$ is the universal object for all quasitoric manifolds determined by combinatorial data $(P,\Lambda)$. Any quasitoric manifold that is determined by a pair $(P,\Lambda)$ is equivariantly homeomorphic to a quotient space of the topological manifold $\tilde \Z_P$ by a free torus action.

Let $F = F_I = \bigcap\limits_{i\in I} F_i$ be a face of codimension $k$ of the polytope $P$. Recall that the embedding $K\hookrightarrow \m T^m$ is given by the $(m\times (m-n))$-matrix ~$C$. Denote by $C_I$ the matrix obtained from $C$ by removing rows with indices from the set $I$.

\begin{lemma}
\cite{BPR07}
\label{dirsum}
The matrix $C_I$ determines a monomorphism of $\m Z^{m-n}$ to a direct summand in $\m Z^{m-|I|}$.
\end{lemma}  

{\it Proof.} 
It is enough to show that the corresponding composition map $K\hookrightarrow \m T^m \to T_{[1,m]\setminus I}$ is a monomorphism. The kernel of projection map $\m T^m\to T_{[1,m]\setminus I}$ is the torus $T_I\subset \m T^m$, and since $K=\ker l$ and $l\colon \m T^m\to \m T^n$ is a characteristic homomorphism, the subgroups $K$ and $T_I$ have trivial intersection. 
$\Box$

\begin{lemma}
Quasitoric manifold $M$ has a canonical $T^n$-invariant Riemannian metric.
\end{lemma}

{\it Proof.}
The action of the compact torus $\m T^m$ on $\m C^m$ is orthogonal, so is
the action of $K\subset \m T^m$ on the moment-angle manifold $\Z_P\subset \m C^m$. The metric on the quotient $M$ = $\Z_P/K$ is therefore well-defined.
$\Box$

\section{Equivariant embeddings of quasitoric manifolds.}

Every vector $a~ = (a_1, \ldots, a_m) \in \m Z^m$ determines a real-algebraic monomial function ${\varphi_{a}\colon \m C^m \to \m C}$ given by the formula:
$$
\varphi_{a}(z_1, \ldots, z_m) = \hat   z_1^{a_1} \cdot \ldots \cdot \hat  z_m^{a_m},
$$
where
\begin{itemize}
{\item $\hat  z_i^{a_{i}} = 1$ if $a_{i} = 0$,}
{\item $\hat z_i^{a_{i}} = z_i^{a_{i}}$ if $a_{i} > 0$,}
{\item $\hat z_i^{a_{i}}  = \bar z_i^{-a_{i}}$ if $a_{i} < 0$.}
\end{itemize}

Let $t=(t_1, \ldots, t_r) \in \m T^r$, $a~= (a_1, \ldots, a_r)\in\m Z^r$, $A$ be an integral $(l\times r)$-matrix and $A_k$ be the $k$-th row of the matrix $A$. We will use the notation $t^{a} = t_1^{a_{1}}\cdot\ldots\cdot t_r^{a_{r}}\in \m T^1$ and $t^A = (t^{A_1}, \ldots,  t^{A_l}) \in \m T^l$. 

It follows directly from the definition that if $a\in \m Z^m$ and $t\in \m T^m$, then $\varphi_a(tz)$ = $t^{a}\varphi_{a}(z)$ for all $z\in\m C^m$. The standard lattice $\m Z^m\subset \m R^m$ may be identified with the lattice of characters of the torus $\m T^m$, which is a subset of conjugate Lie algebra $\mathfrak t_m^*$. The group $K$ is a direct factor in the torus $\m T^m$, so the epimorphism of conjugate Lie algebras $\mathfrak t_m^*\to \mathfrak k^*$ induces an epimorphism of corresponding lattices of characters. Any character $K\to\m T^1$ may be split into a composition $K\hookrightarrow \m T^m\to\m T^1$ and therefore is determined by a vector $a\in \m Z^m$. Conversely, any vector $a\in\m Z^m$ determines a character $k_a\colon K\hookrightarrow \m T^m\to\m T^1$. 

Recall that the embedding of the group $T^{m-n}\simeq K\hookrightarrow\m T^m$ is given by the matrix $C$, so any element $t\in K$ has the form $t = \tau^C$ for some (unique) element $\tau\in\m T^{m-n}$. 

\begin{lemma}
For any $t\in K$ we have $$\varphi_a(t z) = \tau^{a^T C}\varphi_a(z),$$ where $t=\tau^C$ and $\tau\in\m T^{m-n}$.
\end{lemma}
{\it Proof.}
We have $\varphi_a(t z) = t^a\varphi_a(z) = (\tau^C)^a~\varphi_a(z) = \tau^{(a^T C)} \varphi_a(z)$.
$\Box$

\begin{lemma}
\label{samecharacter}
Let $a,a'\in\m Z^m$. Then the following statements are equivalent:
\begin{enumerate}
{\item vectors $a$ and $a'$ generate the same character of the group $K$,}
{\item $C^T a$ = $C^T a'$,}
{\item $(a~- a')\in\op{Im}\Lambda^T$.}
\end{enumerate}
\end{lemma}

{\it Proof.}
The character of the group $K$ induced by the vector $a-a'$ is trivial if and only if $(a-a')\in\ker C^T$. Since $C^T$ is a cokernel matrix for $\Lambda^T$, spaces $\ker C^T$ and $\op{Im}\Lambda^T$ coincide.
$\Box$

Therefore, the standard lattice $\m Z^m\subset \m R^m$ splits into affine layers of the form $(C^T)^{-1}(b)$, $b\in\m Z^{m-n}$, and vectors from the same layer determine the same character $K\to \m T^1$. The sublattice $\ker C^T\cap \m Z^m$ corresponds to the trivial character of $K$. 

\begin{lemma}
\label{vectfamily}
Suppose that vectors $a_1, \ldots, a_q\in \m Z^m$ define the same character $\tilde k\colon K\to \m T^1$. Let $\varphi\colon \Z_P\to \m C^q$ be a restriction of the monomial map $(a_1,\ldots,a_q)$ to the moment-angle manifold $\Z_P\subset \m C^m$. 
\begin{enumerate}
{\item If the character $\tilde k$ is trivial, then $\varphi$ is constant on orbits of the action of $K$ on $\Z_P$. The map $\varphi$ induces a smooth map ${\tilde\varphi\colon M\to \m C^q}$ equivariant with respect to some representation $\m T^n\to \m T^q$, where action of $\m T^q$ on $\m C^q$ is supposed to be standard.}
{\item If the maps $\varphi_{a_1},\ldots, \varphi_{a_q}$ don't vanish simultaneously on $\Z_P$, then the map $\varphi$ induces a smooth map $\tilde\varphi_{\m P}\colon M\to \m CP^{q-1}$ equivariant with respect to some representation $\m T^n\to \m T^q$, where the action of $\m T^q$ on $\m CP^{q-1}$ is supposed to be standard. }
\end{enumerate} 
\end{lemma}

{\it Proof.}
The map $\varphi$ is equivariant with respect to the representation $\m T^m\to \m T^q$ generated by vectors $a_1,\ldots, a_q$.
As mentioned above, $\varphi_a(tz)$ = $\tilde k(t)\varphi_{a}(z)$ for every $t\in K$. So if the character $\tilde k$ is trivial, the map $\varphi\colon \Z_P\to\m C^q$ is constant on orbits of the smooth action of the compact group $K$. The map $\varphi$ therefore induces a smooth map $\tilde\varphi\colon M\to\m C^q$, which is equivariant with respect to some representation $\m T^m/K \simeq \m T^n \to \m T^q$.

If the character $\tilde k$ is nontrivial, the action of any element $t\in K$ on $\Z_P$ multiplies all components of $\varphi$ by the same complex number. Since the image $\varphi(\Z_P)\subset \m C^q$ is compact and does not contain the origin, we may apply the same arguments as above to the composition map $\tilde\varphi_{\m P}\colon\Z_P\stackrel{\varphi}{\longrightarrow} \m C^q\setminus\{0\}\to \m CP^{q-1}$.
$\Box$

\begin{definition}
\label{monom-quasi}
Maps of the form $\tilde\varphi\colon M\to\m C^q$ and $\tilde\varphi_{\m P}\colon M\to\m CP^{q-1}$ constructed by a family of vectors $a_1, \ldots, a_q$ as in lemma \ref{vectfamily} are called {\it monomial maps} of the quasitoric manifold $M$.
\end{definition}

Therefore, functions of the form $\varphi_a$, $a\in \m Z^m$, generate a family of smooth equivariant maps of the manifold $M$. As we will show, this family is wide enough to construct equivariant {\it embeddings} of quasitoric manifolds to projective spaces in terms of combinatorial data $(P,\Lambda)$. These embeddings are parametrized by characters of the group $K$; the trivial character defines an equivariant embedding to Euclidean space.

Consider a generic situation of a compact Lie group $G$ acting smoothly on a compact manifold~$M$. According to Mostow-Palais theorem (\cite{mostow}, \cite{palais}), there exists a smooth embedding $M\to \m R^N$ equivariant with respect to some representation $G\to \m GL(N,\m R)$. The proof uses functional analysis: matrix elements of finite-dimensional representations of $G$ are dense in space of all smooth functions on $M$ with the canonical action of group $G$. Note that Mostow-Palais theorem does not provide any bounds on the dimension of space of the embedding. We will show that for quasitoric manifolds one can construct an equivariant embedding in an explicit way by using monomial maps defined above. A natural bound for the dimension of space of the embedding is in this case given by the number of edges of the polytope $P$. 

Let us fix some character $\tilde k\colon K\to\m T^1$. Since the matrix $C$ determines the isomorphism of $\m T^{m-n}$ and $K$, the character $\tilde k$ is uniquely determined by some vector $\tilde b\in\m Z^{m-n}$. Recall that for a given set of indices $I\subset [1,m]$ we denote by $C_I$ the matrix obtained from $C$ by removing rows with indices from $I$.

\begin{itemize}
{\item Let $r\subset P$ be an edge of the polytope $P$. Since the action of $\m T^n$ on $M$ is locally standard, the isotropy subgroup $\chi(\pi^{-1}(r))$ of the submanifold $\pi^{-1}(r)$ is an ${(n-1)}$-dimensional toric subgroup in $\m T^n$. Denote by $\Phi_r\colon \m T^n\to \m T^1$ the character with kernel isomorphic to $\chi(\pi^{-1}(r))$, and by $\mu_r$  the composition ${\mu_r\colon \m T^m \stackrel{l}{\Arrow} \m T^n \stackrel{\Phi_r}{\Arrow} \m T^1}$. There is an ambiguity in the definition of $\mu_r$ which will be resolved a bit later.
}
{\item Let $\tilde k\colon K\to \m T^1$ be a character of the group $K$ and ${v=F_I\in P}$ be a vertex of $P$. By lemma \ref{dirsum}, the projection map ${p_I\colon K\hookrightarrow \m T^m\to T_{[1,m]\setminus I}}$ is an isomorphism. One may associate with the vertex $v$ a character $k_v\colon \m T^m\to T_{[1,m]\setminus I} \stackrel{p_I^{-1}}{\Arrow} K\stackrel{\tilde k}{\Arrow} \m T^1$.}
{\item Let $\tilde k\colon K\to \m T^1$ be a character of the group $K$, $r\subset P$ an edge of $P$ and $v\in r$ a vertex lying on the edge $r$. Then the pair $(v, r)$ determines a character $k_{v,r}\colon\m T^m \to\m T^1$ defined as the sum of the characters ${k_v\colon \m T^m\to \m T^1}$ and ${\mu_r\colon \m T^m\to \m T^1}$.}
\end{itemize}

The set of pairs of the form $\{$character $\tilde k\colon K\to\m T^1$, vertex $v\in P\}$ determines the set of characters $\{k_v|v\in P\}\subset \op{ch(\m T^m)}$ and the set of triples of the form $\{$character $\tilde k\colon K\to\m T^1$, edge $r\subset P$, vertex $v\in r\}$ determines the set $\{k_{v,r}|v\in r\subset P\}\subset \op{ch(\m T^m)}$, where $\op{ch(\m T^m)}$ is the set of all characters of $\m T^m$. Therefore, every character $\tilde k\colon K\to\m T^1$ determines the set $X_{\tilde k} = (\{k_v\}\cup\{k_{v,r}\}) \subset \op{ch}(\m T^m)$. We list some of the properties of the set $X_{\tilde k}$ and characters $k_v$ and $k_{v,r}$:

\begin{enumerate}
{\item The number of different characters in $X_{\tilde k}$ does not exceed $f_0(P)(n+1)$, because there are at most $f_0(P)$ characters $\{k_v\}$ and at most $n f_0(P) = 2 f_1(P)$ characters of the form $\{k_{v,r}\}$. In practice, $|X_{\tilde k}|$ is often much less than this upper bound, since characters $k_v$ and $k_{v,r}$ may be equal for different vertices $v\in P$ and pairs $(v,r)$.}
{\item The restriction of any of characters $k_v$ and $k_{v,r}$ to the subgroup $K\subset \m T^m$ is equal to $\tilde k\colon K\to \m T^1$.}
{\item For any vertex $v\in P$ the character $\tilde k_v$ is trivial if and only if the character $\tilde k$ is trivial.} 
{\item The character $\mu_r$ does not depend on $\tilde k$ and is nontrivial for every $r\subset P$. The restriction of $\mu_r$ to the subgroup $K$ is trivial.}
{\item Every character $k_v$, $v\in P$, is well-defined, but characters $\Phi_r\colon \m T^n\to\m T^1$ and $\mu_{r}\colon\m T^m\to\m T^1$ are defined only up to multiplication by $\pm1$, so the definition of $k_{v,r}$ is still ambiguous. If $v_0,v_1$ are vertices lying on an edge $r\subset P$, then, as we will show later, the character $k_{v_1} - k_{v_0}$ is a multiple of $\mu_r$. We set $k_{v_0, r} = k_{v_0} + \mu_r$, where $\mu_r$ has the same direction as $k_{v_1} - k_{v_0}$. If $k_{v_1} = k_{v_0}$, then we assume that the first nonzero coordinate of $m$-vector defining $\mu_r\colon\m T^m\to \m T^1$ is positive.}
{\item If the character $\tilde k$ is trivial, the cardinality of $X_{\tilde k}$ does not exceed $f_1(P) + 1$, because $k_v\equiv 1$ for all $v\in P$ and if $v_0,v_1$ are vertices of an edge $r\subset P$, then $k_{v_0, r} = k_{v_1,r} = \mu_r$.}
\end{enumerate}

Vectors $w_r$ that define characters $\Phi_r$, $r\subset P$, form an integral $(n\times q)$-matrix $W$. Consider the corresponding linear representation $\Phi$ of the torus $\m T^n$ in linear space $\m C^q$. The moment map $\pi\colon M\to P\subset \m R^n$ is equivariant with respect to the trivial torus action on $\m R^n$. 

\begin{theorem}
\label{main}
The moment map $\pi\colon M\to P$ can be extended to $\mbox{a~real-algebraic}$ embedding $\pi\times\tilde\varphi\colon M \to \m R^n\times \m C^q$ equivariant with respect to the representation $\Phi\colon\m T^n\to\m T^q$.
\end{theorem}

The number $q$ in theorem \ref{main} does not exceed the number of edges of $P$, as follows from the construction of the representation $\Phi$. 

Theorem \ref{main} is a corollary of a more general result on projective embeddings of a quasitoric manifold $M$. Every character from the set $X_{\tilde k}$ constructed above determines a vector $a\in\m Z^m$, which in turn determines a monomial function $\varphi_a\colon \m C^m\to \m C$. Therefore, the set $X_{\tilde k}$ determines a monomial map $\varphi_{\tilde k}\colon \m C^m\to \m C^q$, where $q = |X_{\tilde k}|$. As we will show, the corresponding map $\tilde\varphi_{\m P, \tilde k}\colon M\to \m CP^{q-1}$ is well-defined. 

Every representation $\tilde k\colon K\to\m T^1$ determines a complex linear bundle $\xi_{\tilde k}\to M$, where ${\xi_{\tilde k} = \Z_P\times_K \m C}$ and the action of $K$ on $\m C$ is defined via the character $\tilde k$. It is known from algebraic topology that the groups $H^2(M,\m Z)$ and $\op{ch}(K)$ are isomorphic and so we may assign to the character $\tilde k\colon K\to \m T^1$ the cohomology class $c_1(\xi_{\tilde k})\in H^2(M,\m Z)$. Note that the class $c_1(\xi_{\tilde k})$ is determined by a map $f\colon M\to\m CP^N$ such that the bundle $\xi_{\tilde k}$ over $M$ is a pull-back of the tautological linear bundle on $\m CP^N$ for $N$ large enough. It turns out that the map $f$ may be extended to an embedding of $M$.

\begin{theorem}
\label{projmain}
The set of characters $X_{\tilde k}$, where $\tilde k\colon K\to\m T^1$ is an arbitrary character, determines a monomial map $\tilde\varphi_{\m P, \tilde k}\colon M\to\m CP^{q-1}$ (see def. \ref{monom-quasi}) that can be extended to an~embedding ${\pi\times\tilde\varphi_{\m P, \tilde k}\colon M \to P\times \m CP^{q-1}}$. The induced cohomology pullback ${\tilde\varphi_{\m P, \tilde k}^*\colon H^2(\m CP^{q-1}, \m Z)\to H^2(M, \m Z)}$ coincides with the classifying map $H^2(\m CP^{\infty}, \m Z) \stackrel{\simeq}{\Arrow} H^2(\m CP^{q-1}, \m Z) \to H^2(M, \m Z)$ of the bundle $\xi_{\tilde k}$. If the character $\tilde k$ is trivial, then the image of the map $\tilde\varphi_{\m P,\tilde k}$ lies entirely in some affine chart $\m C^{q-1}\subset\m CP^{q-1}$.
\end{theorem}

\begin{remark}
If the pair $(P,\Lambda)$ determines a projective toric variety, then embeddings in theorem \ref{projmain} generalize projective embeddings from algebraic geometry.
\end{remark}

Let us provide an explicit algorithm for constructing the set ${X_{\tilde k}\subset\op{ch}\m T^m}$ starting from $(P,\Lambda)$ and a character $\tilde k\colon K\to\m T^1$.

\begin{enumerate}
{\item Construct the matrix $C$ that fits into the short exact sequence $$1\to\m T^{m-n} \stackrel{C}{\Arrow} \m T^m \stackrel{\Lambda}{\Arrow} \m T^n \to 1.$$ Recall that $\op{Im}(C) = K\subset \m T^m$.}
{\item Using $C$, we may identify the character $\tilde k\colon K\to\m T^1$ with a~character of standard torus $\m T^{m-n}$ given by a vector $\tilde b\in \m Z^{m-n}$.}
{\item For a given vertex $v = F_I\in P$ consider an $((m-n)\times(m-n))$-matrix $C_I$ obtained from $(m\times(m-n))$-matrix $C$ by removing $n$ rows corresponding to the index set $I$ (by lemma \ref{dirsum}, the matrix $C_I$ is invertible).} 
{\item Consider a vector $b_v\in \m Z^m$ obtained from the vector ${(C^T_I)^{-1}(\tilde b)\in \m Z^{m-n}}$ by putting zeros in places from the index set $I$.  Let $k_v\colon\m T^m\to\m T^1$ be a character corresponding to the vector $b_v$.}
{\item For every edge $r\subset P$ and vertex $v\in r$ denote by $v'\in P$ an opposite vertex on that edge. Suppose first that $b_{v'}\ne b_v$. We define vector $a_{v, r}$ as the closest lattice point to $b_v$ lying on the ray starting in $b_v$ and containing $b_{v'}$. If $b_v = b_{v'}$, we set $a_{v, r} = b_v + u_r$, where $u_r = \alpha_r^T \Lambda$ and $\alpha_r$ is a primitive $n$-vector orthogonal to all columns of $\Lambda$ with indices from the index set $J$, $r=F_J$. We require first nonzero component of the vector $u_r$ to be positive. 
}
{\item Every vector $a_{v, r}$, $v\in r\subset P$, defines a character $k_{v,r}\colon\m T^m\to\m T^1$.}
{\item Define the set $X_{\tilde k}$ as the set formed by all characters $k_v$, $v\in P$, and $k_{v,r}$, $v\in r\subset P$.}
\end{enumerate}

Note that if $r\subset P$ is an edge containing vertices $v$ and $v'$, the vector $u_r = \alpha_r^T \Lambda$ is collinear with $b_{v'} - b_v$, if the last vector is not zero.


\begin{example}
Consider a manifold $\m CP^n$ that is a toric variety over the polytope ${\Delta^n\subset \m R^{n+1} = \m R^m}$. In this case $(n\times (n+1))$-matrix $\Lambda$ has the form
$$
\Lambda
=
\begin{pmatrix}
1 & \ldots & 0 & -1\\
\vdots & \ddots & \vdots & \vdots\\
0 & \ldots & 1 & -1\\
\end{pmatrix},
$$
and the matrix $C$ is a vector $(1,\ldots, 1)^T$ of length $(n+1)$. The characters of the group $K$ are defined by numbers $\tilde b\in \m Z^{m-n}=\m Z$.

Let us first consider the case when $\tilde k$ is a trivial character and $\tilde b = 0$. For every vertex $v=F_I\in \Delta^n$ we have $C_I=(1)$. Then ${b_v = (0,\ldots,0)\in\m Z^{n+1}}$ for all vertices $v\in\Delta^n$. If $r=F_J\subset \Delta^n$ is an edge, then vector $\alpha_r\in\m Z^n$ may be found explicitly, by solving $(n-1)$ linear equations. Suppose, for example, that $r=F_J$, where $J=\{1,\ldots, n-1\}$. Then the primitive vector $\alpha_r$ is orthogonal to first $(n-1)$ columns of the matrix $\Lambda$, which coincide with first $(n-1)$ columns of unit $(n\times n)$-matrix.
It follows that $\alpha_r=(0,\ldots,0,1)\in\m Z^n$. The vector $u_r = \alpha_r^T\Lambda\in\m Z^{n+1}$ is equal to $(0,\ldots, 0, 1, -1)$. 

Similarly, $u_r = \alpha_r^T\Lambda= (e_i - e_j)$ for every edge ${r=F_J}$, ${J = [1,n+1]\setminus\{i,j\}}$, ${1\leqslant i < j \leqslant n+1}$. So ${X_{\tilde k} = \{\op{c}_0\} \cup \{\op{c}_{e_i-e_j} \mid 1\leqslant i < j \leqslant n+1\}}$, where we denote by $\op{c}_x\in\op{ch}(\m T^m)$ the character determined by vector $x\in\m Z^m$.  In general, for every projective toric variety the vector $u_r$, $r\subset P$, coincides with primitive vector collinear with the edge $i_P(r)\subset \m R^m$.  

We have $\varphi_{e_i-e_j}=z_i\bar z_j$ for $1\leqslant i<j \leqslant n+1$. It is clear from the construction that maps $\varphi_{e_i-e_j}$ are invariant under the diagonal circle action on $\m C^{n+1}$. The moment-angle manifold $\Z_{\Delta^n}$ of the simplex $\Delta^n$ is the sphere $S^{2n+1}$. The projection map $\rho_{\Delta^n}$ is given by the formula $(z_0,\ldots, z_n)\to (z_0\bar z_0, \ldots, z_n\bar z_n)$. We may identify space $\m R^{n+1}\times \m C^{\frac{n(n+1)}2}$ with space of Hermitian $((n+1)\times(n+1))$-matrices. Then the map
 $\rho_{\Delta^n}\times\varphi\colon S^{2n+1}\to \m R^{n+1}\times\m C^{\frac{n(n+1)}2}$ is given by the formula
$$
z \mapsto z^T \bar z,
$$
where $z = (z_0, \ldots, z_n)$, $\bar z = (\bar z_0, \ldots, \bar z_n)$ and $z^T$ is $z$ transposed. 

The map $\rho_{\Delta^n}\times\varphi\colon S^{2n+1}\to \m R^{n+1}\times\m C^{\frac{n(n+1)}2}$ generates an embedding of the quasitoric manifold $${\pi\times\tilde\varphi\colon\m CP^n\to \m R^{n+1}\times\m C^{\frac{n(n+1)}2}}.$$ Note that the total number of monomial functions is equal to the number of edges in $\Delta^n$. This number is less than the upper bound for $|X_{\tilde k}|$, which is $f_0(\Delta^n)(n+1)$ = $(n+1)^2$.

Consider now the simplest non-trivial character of the subgroup $K$ determined by the element $\tilde b = (1)\in \m Z$. Since $C_I=(1)$ for any vertex $v=F_I\in P$, the character $k_v$ corresponding to ${v = e_i = (0, \ldots, 0, 1, 0, \ldots, 0)\in \Delta^n}$ is also determined by the vector $e_i$, $i\in [1,n+1]$. So all of the vectors $b_v, v\in P$, are pairwise distinct. If $v' = e_j$ is another vertex of the simplex, then the edge $r\subset \Delta^n$ containing $v$ and $v'$ does not contain any other lattice points, so $a_{v, r} = b_{v'}$.

Therefore, the map $\tilde\varphi_{\m P,\tilde k}$ is formed by $(n+1)$ monomial functions which in this case coincide with coordinate functions on $\m C^{n+1}$. The corresponding projective map is simply the identical map $\m CP^n\to\m CP^n$. We can see that the dimension of the projective embedding is even less than the dimension of the affine embedding constructed above. Another important thing to note is that the embedding is generated only by the component  $\tilde\varphi_{\m P, \tilde k}$: the moment map $\pi$ is not necessary in this example. This is common for quasitoric manifolds coming from projective toric varieties; we will give the precise formulation later, in \ref{projtoric}.
\end{example}

In the remaining part of this section we prove theorem \ref{projmain}. Recall that for any vertex $v\in P$ $b_v\in\m Z^m$ is the vector corresponding to the character $k_v\in X_{\tilde k}$ and $a_{v,r}\in\m Z^m$ the vector corresponding to the character $k_{v,r}\in X_{\tilde k}$. 
We show first that the projective map $\tilde\varphi_{\m P, \tilde k}$ is $\mbox{well-defined}$. It is enough to check that there exists a subset in $X_{\tilde k}$ such that the corresponding monomial functions never vanish simultaneously on the moment-angle manifold $\Z_P\subset \m C^m$.

\begin{lemma}
\label{vertnonzero}
The maps $\varphi_{b_v}$, $v\in P$, do not vanish simultaneously on $\Z_P$.
\end{lemma}

{\it Proof.}
If $\tilde k\equiv 1$, then $\varphi_{b_v}\equiv 1$ for all $v\in P$, so the statement of the lemma is obviously true. Suppose that the character $\tilde k$ is non-trivial and all maps $\varphi_{b_v}$ vanish in some point $z\in \Z_P$. Consider now some vertex $v_0\in F(\rho_P(z))$. If $F(\rho_P(z))$ = $F_I$, then all coordinates of the vector $b_{v_0}$ with indices from the set $I$ are zero. Therefore, the monomial function $\varphi_{b_{v_0}}$ depends only on nonzero coordinates of the point $z$.
$\Box$

We see that the map $\tilde\varphi_{\m P, \tilde k}\colon M\to \m CP^{r-1}$ is well-defined. Suppose now that $(\pi\times\varphi)(z) =  (\pi\times\varphi)(z')$ for some $z, z'\in \Z_P$. Since ${\pi(z) = \pi(z')}$, the points $z$ and $z'$ lie in the same orbit of the action of $\m T^m$, so there exists an element ${t = (t_1, \ldots, t_m)}\in\m T^m$ such that $z'_i = t_i z_i$, $i\in [1, m]$. If ${F(\rho_P(z)) = F(\rho_P(z'))}$ = $F_I$, then this element $t\in \m T^m$ is defined only up to the coordinates from the set $I$, which, in turn, may be arbitrary, since $z'_i$ = $z_i$ = $0$ for $i\in I$.
Showing that $\pi\times\varphi_{\m P, \tilde k}$ is an embedding of the quasitoric manifold $M$ is equivalent to proving that if $\pi(z)=\pi(z')$ and $\tilde\varphi_{\m P,\tilde k}(z)=\tilde\varphi_{\m P, \tilde k}(z')$, then there exists an element $t\in K$ such that $z'$ = $tz$.

Choose an arbitrary vertex $v_0\in F_I$ and let $r_1,\ldots, r_{\dim F_I}$ = $r_{n-|I|}$ be edges of the face $F_I$ containing $v_0$.

\begin{lemma}
\label{different}
If $\dim F_I >1$, then isotropy subgroups $\chi(r_1), \ldots, \chi(r_{\dim F_I})$ are pairwise distinct. 
\end{lemma}

{\it Proof.}
The statement of the lemma follows from the fact that the action of $\m T^n$ is locally standard in the neighborhood of the point $\pi^{-1}(v_0)\in M$.
$\Box$

It follows that vectors $u_{r_i}$ = $a_{v_0,r_i}-b_{v_0}$, $i\in [1, \dim F_I]$, are also pairwise distinct and linearly independent. 

\begin{lemma}
\label{unitmatrix}
If $v_0$ = $\bigcap\limits_{j=1}^{n} F_j$, then $j$-th coordinate of the vector $u_{r_i}$ is equal to $\pm1$, if $v_0$ = $r_i\cap F_j$, and to  zero otherwise (if $r_i\subset F_j$).
\end{lemma}

\begin{corollary}
\label{coordzero}
All of the coordinates of $u_{r_i}$, $i\in [1, \dim F_I]$, corresponding to the indices from the set $I$ are zero. 
\end{corollary}

{\it Proof.}
As we know, every vector $u_{r_i}$ determines a character $\mu_{r_i}\colon \m T^m\to \m T^1$ such that ${t\in \ker\mu_{r_i}}$ if and only if $l(t)$ acts trivially on the two-dimensional sphere $\pi^{-1}(r_i)$. If we have $r_i\subset F_j$, then $j$-th coordinate of all points in $\rho_P^{-1}(r_i)\subset \Z_P\subset\m C^m$ is zero. So $s\in\ker\mu_{r_i}$ for every $s\in T_j$. Since the homomorphism $l\colon \m T^m\to \m T^n$ is a characteristic function, the restriction of $l$ to the torus generated by subgroups $T_j$, $v_0\in F_j$, is an isomorphism. If $v_0$ = $r_i\cap F_j$ and $j$-th coordinate of $u_{r_i}$ is not $\pm 1$, then the isotropy subgroup of any point in $\pi^{-1}(\stackrel{\circ}{r_i})$ cannot be isomorphic to $\m T^{n-1}$, so the action of $\m T^n$ on $M$ would not be locally standard in this case. $\Box$

If $z'$ = $tz$ and $\tilde\varphi_{\m P, \tilde k}(z')$ = $\tilde\varphi_{\m P, \tilde k}(z)$, then $\varphi_{u_{r_i}}(z')$ = $\varphi_{u_{r_i}}(z)$ for all $u_{r_i}$. Denote by $p_I(t)$ an image of the element ${t\in \m T^m}$ under the canonical projection $p_I\colon \m T^m\to T_{[1,m]\setminus I}$. Then ${\varphi_{u_{r_i}}(z') = \varphi_{u_{r_i}}(tz) = p_I(t)^{u_{r_i}}\varphi_{u_{r_i}}(z)}$ and we have $p_I(t)^{u_{r_i}}$ = $1$ for all $i \in [1,\dim F_I]$. We will show that $p_I(t)\in p_I(K)$; surely, this implies that there exists an element $t\in K$ such that $z' = tz$.

The vectors $u_{r_i}^T$, $i=1\ldots \dim F_I$, with coordinates from the index set $I$ removed form an integral matrix $\Lambda_I$ of the size $(\dim F_I\times (m-|I|))$.

\begin{lemma}
\label{exact}
The sequence
$$
1 \to \m T^{m-n} \stackrel{C_I}{\longrightarrow} T_{[1,m]\setminus I} \stackrel{\Lambda_I}{\longrightarrow} \m T^{\dim F_I} \to 1
$$
is exact.
\end{lemma}

{\it Proof.}
By \ref{coordzero}, we have $\Lambda_I C_I$ = $0$, since $u_{r_i}^T C$ = $0$ by lemma \ref{samecharacter}. By lemma \ref{dirsum}, the matrix $C_I$ determines a monomorphism of the torus $\m T^{m-n}$ to a direct factor in $T_{[1,m]\setminus I}$. The matrix formed by $u_{r_i}$ contains a~unit submatrix of rank $(n-|I|)$ by the lemma \ref{unitmatrix}. Therefore, the map $T_{[1,m]\setminus I} \stackrel{\Lambda_I}{\longrightarrow} \m T^{\dim F_I}$ is epimorphic and its kernel coincides with $\op{Im}(C_I)$.
$\Box$

Lemma \ref{exact} completes the proof of the part of theorem \ref{projmain} concerning the injectivity of $\pi\times\tilde\varphi_{\m P, \tilde k}$. Indeed, if $p_I(t)^{u_{r_i}}$ = $1$ for all ${r_i}$, then $p_I(t)\in C_I(\m T^{m-n})$ = $p_I(K)$. 

Let us now prove the part of theorem \ref{projmain} concerning pullback map in cohomology. If ${\m CP^{q-2}\subset \m CP^{q-1}}$ is an embedding of one of standard coordinate hyperplanes, then the corresponding monomial function is given by a vector ${a\in \m Z^m}$, where $a$ is either one of the vectors $b_v$, $v\in P$, or one of the vectors $a_{v,r}$, $v\in r\subset P$. The preimage $\tilde\varphi_{\m P, \tilde k}^{-1}(\m CP^{q-2})\subset M$ is then formed by a union of characteristic submanifolds ${\pi^{-1}(F_j)\subset M}$, $F_j\subset P$. The multiplicity of each component $\pi^{-1}(F_j)$ is equal to $j\mbox{-th}$ coordinate of the vector $a$. The second cohomology group of the quasitoric manifold $M$ is a free abelian group generated by duals of characteristic submanifolds and factorized by relations given by rows of the matrix $\Lambda$. If the character $\tilde k$ is non-trivial, then, by lemma \ref{samecharacter}, $a\notin\op{Im}(\Lambda^T)$, so the pullback map is nonzero. If the chacacter $\tilde k$ is trivial, then the image of the map $\tilde\varphi_{\m P, \tilde k}(M)$ lies entirely in an affine chart of $\m CP^{q-1}$.

Let us finally check that the map $\pi\times\tilde\varphi_{\m P, \tilde k}$ is a smooth embedding of $M$. The map $\rho_P\times\varphi_{\m P, \tilde k}$ is constant on orbits of the action of $K$, so any vector that is tangent to orbit lies in the kernel of $(\rho_P\times\varphi_{\m P, \tilde k})_*$. We need to show that the kernel of $(\rho_P\times\varphi_{\m P, \tilde k})_*$ in any point $z\in\Z_P$ does not have vectors other than tangent to orbit. The map $\rho_P\times\varphi_{\m P, \tilde k}$ commutes with the action of $\m T^m$ on $\Z_P$, so it is enough to check the case $z\in s(P)\subset \Z_P$, i.e. all coodinates of $z$ are real and non-negative. Let $v$ be an arbitrary vertex of the face $F = F_I = F(\rho_P(z))$. We may suppose that $I=\{1\ldots k\}$ and $v = F_1\cap\ldots\cap F_n$.

Consider the following system of local coordinates in the point ${z\in\Z_P}$:
\begin{enumerate}
{\item Cartesian coordinates $x_1, y_1, \ldots, x_k, y_k$;}
{\item arguments $t_{k+1}, \ldots, t_n$ of complex coordinates $z_{k+1}, \ldots, z_n$;}
{\item local coordinates in the image $s(\stackrel{\circ}{F})\subset\Z_P$.}
\end{enumerate}

Choosing a system of local coordinates in the orbit of the action of $K$ gets us a complete coordinate system in the neighborhood of $z\in\Z_P$.

Let $r_1, \ldots, r_n\subset P$ be edges outgoing from $v$. From the proof of the lemma \ref{vertnonzero} it follows that $\varphi_{b_v}(z)\ne 0$. Now consider the following local coordinates in the image of the map $\rho_P\times\varphi_{\m P, \tilde k}$:
\begin{enumerate}
{\item real and imaginary parts of monomial functions $\varphi_{a_{v, r_i}} / \varphi_{b_v}$, ${i=1\ldots k}$;}
{\item arguments of monomial functions $\varphi_{a_{v, r_i}} / \varphi_{b_v}$, $i=k+1,\ldots, n$;}
{\item local coordinates in the neighborhood of the point ${\rho_P(z) = \pi(\pi_l(z))}$.}
\end{enumerate}

By lemma \ref{unitmatrix}, every function $\varphi_{a_{v, r_i}} / \varphi_{b_v}$, $i=1\ldots n$, has the form $\hat z_i^{\pm1} F(z)$, where $F$ is a monomial function depending only on variables $z_{n+1}^{\pm1}, \bar z_{n+1}^{\pm1},\ldots, z_m^{\pm1}, \bar z_m^{\pm1}$. This implies that the square matrix of all partial derivatives in provided coordinate systems is non-degenerate.

\section{Projective embeddings of toric varieties}

It is well known that a projective toric variety $M$ is a symplectic quotient of a corresponding moment-angle manifold $\Z_P\subset \m C^m$, which coincides with the level set of a Hamiltonian $\m C^m\to\m R^{m-n}$ defined by the action of $K\subset \m T^m$ on $\m C^m$.
By Delzant theorem (\cite{delzant}), the vertices of polytope $P$ lie in the lattice $\m Z^n\subset\m R^n$.  The lattice points from the interior of $P$ (including those on the boundary) generate an~embedding of $M$ to the projective space indexed by these points. In other words, the statement of the theorem \ref{projmain} may be strenghtened for toric varieties -- the embedding is determined solely by monomial functions, without the moment map. In this section we will give the proof of this statement which will not use algebraic geometry. This allows to construct projective embeddings of symplectic toric varieties using only their symplectic quotient construction. 

By Kodaira theorem (\cite{kodaira}), if $M$ is a compact complex manifold with a positive holomorphic linear bundle (for example, this is true, when $M$ has a rational K\"ahler form), then there exists a complex-analytic embedding of $M$ to $\m CP^n$. A generalization of this theorem to the class of compact symplectic manifolds is Gromov-Tishler theorem (\cite{gromov}, \cite{tischler}): any symplectic manifold with a rational symplectic form $\omega$ may be embedded into $\m CP^n$ in such a way that the form $\omega$ will be a pullback of the standard symplectic form on $\m CP^n$. The proof of both theorems is essentially based on analysis, in particular, on theory of elliptic operators. In the case of quasitoric manifolds the usage of combinatorial data $(P,\Lambda)$ allows to avoid this machinery.

If~$C_P$ is a cokernel matrix for $A_P$ (see (\ref{eq2})), then one may choose $C_P^T$ as the matrix $C$ determining the subgroup ${K\hookrightarrow \m T^m}$. The polytope $P\subset \m R^n$ is mapped under $i_P\colon x\mapsto A_P x + b_P$ into an intersection of the corresponding affine $n$-plane with $\m R^m_{\geqslant 0}$.
The lattice points from the interior of $P$ map into the set of lattice points $y_j\in\m Z^m$ satisfying two conditions:
\begin{itemize}
{\item coordinates of all points $y_j$ are non-negative,}
{\item $C_P(y_i)$ = $C_P(b_P)$.} 
\end{itemize}

This implies that all degrees of all variables $z_i$ of a monomial map generated by any of $y_j\in S_P$ are non-negative. It is a known result in toric geometry that the map defined by all of $y_j\in S_P$ generates an~embedding of $M$ into space $\m CP^{|S_P|-1}$ (e.g.\cite{fulton}). 

We will give the proof of this result in a slightly strengthened form. Since $C^T$ = $C_P$, the points $\{y_j\}$ all determine the same character $\tilde k_P\colon K\to \m T^1$, which is also a restriction of the character defined by $b_P\in \m Z^m$ to the subgroup $K$. Let us construct the character set $X_{\tilde k_P}$ and system of vectors $\{b_v, a_{v,r}\}$ as in theorem \ref{projmain}. In the case of toric variety vectors $b_v$, $v\in P$, correspond to the vertices of $P$, that is, $b_v$ = $i_P(v)\in\m Z^m$ for all $v$. For every vertex $v\in P$ and an edge $r\subset P$ containing $v$ the vector $a_{v,r}$ is simply the point on the edge $i_P(r)$ that is closest to the vertex $i_P(v)$ (and not equal to it). 

Note that all of the vertex vectors $b_v$ are pairwise different, but this might not be true for $a_{v,r}$. If an edge $r$ contains vertices $v_0, v_1$ and there is no lattice points on $i_P(r)$ except $i_P(v_0)$ and $i_P(v_1)$ themselves, then $a_{v_0,r}$ = $b_{v_1}$ and $a_{v_1,r}$ = $b_{v_0}$, as in the example of $\m CP^n$ that we've seen before. If there is exactly one lattice point on $i_P(r)$ except $i_P(v_0)$ and $i_P(v_1)$, then it is the middle of the edge $i_P(r)$ so we have ${a_{v_0,r} = a_{v_1,r}}$. There is also an important difference from the general situation of quasitoric manifold: exactly $n$ coordinates of each vector $b_v$, $v\in P$, are equal to zero and the rest of $(m-n)$ coordinates must be positive. 

We now formulate the main result of this section.

\begin{theorem}
\label{projtoric}
Suppose that combinatorial data $(P,\Lambda)$ of a quasitoric manifold $M$ satisfies following conditions:
\begin{itemize}
{\item $b_P\in\m Z^m$,}
{\item we have $A_P^T = B\Lambda D$, where $B$ is an invertible integral ${(n\times n)}$-matrix 
and $D$ is a diagonal $(m\times m)$-matrix with all nonzero elements equal to $\pm 1$.
}
\end{itemize} 
Then the vector $b_P\in \m Z^m$ determines a character $\m T^m\to\m T^1$ whose restriction to the subgroup $K\subset\m T^m$ gives a character $\tilde k_P\colon K\to \m T^1$. 
The projective map ${\tilde\varphi_{\m P, \tilde k_P}\colon M\to \m CP^{q-1}}$ constructed from the set $X_{\tilde k_P}$ (see \ref{projmain}) is a smooth embedding.
\end{theorem}

\begin{corollary}
Suppose that $M$ is a complete nonsingular toric variety determined by a normal fan of some simple lattice polytope $P$. Consider monomials determined by vectors ${b_v = i_P(v)\in\m Z^m}$, $v\in P$, and $a_{v,r}$, ${v\in r\subset P}$, where $a_{v,r}$  is a closest point to $i_P(v)$ on the edge $i_P(r)$. Denote by $X_P$ the set of all such monomials. Then the corresponding map of $\m C^m$ to $\m C^q$, $q=|X_P|$, induces an equivariant embedding $\tilde\varphi_{\m P, \tilde k_P}\colon M\to \m CP^{q-1}$.
\end{corollary}

It is well known in algebraic geometry that a complete toric variety is projective if and only if its fan is a normal fan of some simple polytope (see e.g.\cite{fulton}). As shown in \cite{suyama}, there exist complete toric varieties that are not quasitoric manifolds.

\begin{example}
Consider a quasitoric manifold over the polytope ${P =\Delta^n}$. In this case we have $m = n+1$, $\Z_P = S^{2n+1}$ and $K$ is a one-dimensional subgroup in $\m T^{n+1}$ determined by $((n+1)\times 1)$-matrix $C$. The condition of $K$ acting freely on $S^{2n+1}$ is equivalent to the fact that $C = (c_1, c_2, \ldots, c_{n+1})^T$, where $c_i=\pm1$. Then the reduced form of $\Lambda$ is:
$$
\Lambda
=
\begin{pmatrix}
1 & \ldots & 0 & -c_1 c_{n+1}\\
\vdots & \ddots & \vdots & \vdots\\
0 & \ldots & 1 & -c_n c_{n+1}
\end{pmatrix}.
$$
Let $B$ be a diagonal $(n\times n)$-matrix with diagonal values $c_1c_{n+1},\ldots, c_n c_{n+1}$, and $D$ be a diagonal $(n+1)\times (n+1)$-matrix with diagonal values $c_1c_{n+1},\ldots, c_n c_{n+1}, 1$. Then we have
$$
B \Lambda D
=
\begin{pmatrix}
1 & \ldots & 0 & -1\\
\vdots & \ddots & \vdots & \vdots\\
0 & \ldots & 1 & -1
\end{pmatrix}
=
A_P^T.
$$ 
Consider the character of $\m T^m$ given by the vector $b_P = (0, \ldots, 0, 1)^T$. The restriction of this character to the subgroup $K\subset \m T^m$ is defined by the single number $C^T(b_P) = c_{n+1}\in \m Z$. Every vertex $v=F_I\in\Delta^n$ defines the $(1\times 1)$-matrix $C_I$ equal to $(c_i)$. Consequently, the vector $b_v$ has the form $(0, \ldots, c_i c_{n+1}, \ldots, 0)\in\m Z^{n+1}$. 

Note that all vectors $b_v$ are pairwise different. If $r$ is an edge containing vertices $v_0$ and $v_1$, then $a_{v_0, r}$ is the closest lattice point to $b_{v_0}$ on the ray that starts in $b_{v_0}$ and goes through $b_{v_1}$. In our case we always have $a_{v_0, r} = b_{v_1}$. If $v_i\in \Delta^n$, $i \in [1, n+1]$, then $\varphi_{b_{v_i}} = (\hat z_i)^{c_i c_{n+1}}$. These monomial functions generate an embedding of quasitoric manifold $M = \Z_{\Delta^n}/K$ to the projective space $\m CP^n$ of the same dimension. Clearly this is an equivariant diffeomorphism with respect to the representation $\m T^{n+1}/K \to \m T^n$ and standard action of $\m T^n$ на $\m CP^n$.

\end{example}

{\it Proof of theorem \ref{projtoric}.}
After a suitable basis change in $n$-dimensional torus $\m T^m / K$ we may assume that $B$ is a unit matrix and $A_P^T = \Lambda D$. Next, since $0 = A_P^T C_P^T = \Lambda D C_P^T$, the matrix $C = DC_P^T$ is a cokernel matrix for $\Lambda$. 

Suppose that the theorem is proved when $D$ is a unit matrix and $\Lambda = A_P^T$. Then the set of vectors $b_v$ and $a_{v,r}$ constructed from the character defined by vector $C_P(b_P)\in\m Z^{m-n}$ determines a monomial map ${\Z_P/\tilde K = \tilde M\to\m CP^{q-1}}$. Here the subgroup $\tilde K\subset\m T^m$ is an image of the torus $\m T^{m-n}$ under the embedding induced by the matrix ${DC = C_P^T}$. It follows that the system of vectors $Db_v$ and $Da_{v,r}$ generates an embedding $M\to\m CP^{q-1}$ of initial quasitoric manifold $M$. Therefore, it suffices to prove theorem \ref{projtoric} in the case of $B$ and $D$ being unit matrices.
Moreover, the proof of theorem \ref{projmain} shows that we only need to establish the following result.

\begin{lemma}
\label{modules}
If $z', z\in \Z_P$ and $\tilde\varphi_{\m P, \tilde k_P}(z')$ = $\tilde\varphi_{\m P, \tilde k_P}(z)$, then $|z_j'|$ = $|z_j|$ for all $j \in [1,m]$. 
\end{lemma}

{\it Proof.}
Suppose that $z',z\in \Z_P$ and $\tilde\varphi_{\m P, \tilde k_P}(z')$ = $\tilde\varphi_{\m P, \tilde k_P}(z)$. Let us first show that there exists a face $F\subset P$ such that $\rho_P(z'), \rho_P(z)\in\stackrel{\circ}{F}$. 

The sets $I'$ and $I$ of indices of zero coordinates of $z$ and $z'$ correspond to some faces $F_{I'}$ and $F_{I}$ of the polytope $P$. If $F_{I'}\ne F_{I}$, then there exists a vertex $v$ that belongs only to one of the faces; suppose that $v\in F_{I'}$, $v\notin F_{I}$ and $F_j$ is a hypersurface containing $F_{I}$ and not containing $v$. Then the degree of $z_j$ in the monomial $\varphi_{b_v}$ is strictly positive, so $\varphi_{b_v}(z)$ = $0$. Since the set of  nonzero coordinates of the vector $b_v$ is contained in the set $[1,m]\setminus I'$, we have $\varphi_{b_v}(z')\ne~0$. $\Box$

As a consequence, there exists an index set $I$ and a face $F_I\subset P$ such that $z'_i = z_i = 0 \Longleftrightarrow i\in I$. Let $J = [1, m]\setminus I$. Consider vectors $\beta$, $\beta'$, $w$, $w'$ with components $\beta_j = |z_j|^2$, $\beta'_j = |z'_j|^2$, $w_j = (\log |z_j|^2)$,  $w'_j = (\log |z'_j|^2)$, $j~\in~J$. Consider the composition of the affine maps ${\m R^k\stackrel{j_F}{\Arrow} \m R^n \stackrel{i_P}{\Arrow} \m R^m}$, where the first map is an affine embedding of ${F\subset\m R^k}$ onto the face $F_I\subset P\subset\m R^n$. Then we have $\beta',\beta \in j_F(F)$.
Let $x',x\in F$ be the points in the interior of $F\subset\m R^k$ such that $i_P(j_F(x')) = \beta'$, $i_P(j_F(x)) = \beta$.

\begin{proposition}
\label{jsympos}
Let $H\colon U\to \m R^k$ be a smooth map defined on an open convex region $U\subset \m R^k$. Suppose that Jacobian $J_H$ of $H$ is symmetric and positive definite in each point of $U$. Then $H$ is an embedding.
\end{proposition}
{\it Proof.} Let $x,y \in U$, $x\ne y$ and $p(t)$ = $x + t(y-x)$ be the segment joining $x$ and $y$. Since $J_H$ is symmetric and positive definite, we have $(J_H(w), w)>0$ for any vector $w$ in any point of $U$. Denote by $H_v$ the projection of $H$ to the direction $v = (y-x) = (v_1,\ldots,v_k)$. Then $H_v = \sum\limits_i v_i H_i$, where $H=(H_1,\ldots,H_k)$. The derivative of $H_v$ in the direction $v$ is therefore equal to $\sum\limits_i v_i \frac{H_i}{\partial v} = \sum\limits_i\sum\limits_j v_i v_j \frac{\partial H_i}{\partial x_j}= (J_H(v), v)$. Since $(J_H(v), v)>0$ everywhere on the path $p(t)$, $H_v(y) > H_v(x)$, so $H(y)\ne H(x)$. $\Box$

\begin{remark}
$J_H$ being symmetrical is equivalent to the fact that the differential form $\omega = \sum\limits_{i=1}^n H_i dx_i$ is closed. By Poincar\'e lemma, the form $\omega$ is exact on $U$. Therefore $\omega = dF$ for some smooth $F\colon U\to \m R$ and in this case $J_H$ is Hessian of $F$. 
\end{remark}

We will show below that the components of the monomial map $\varphi_{\tilde k_P}$ generate the map
${H\colon \stackrel{\circ}{F}\to\m R^k}$ defined on the open face $\stackrel{\circ}{F} \subset \m R^k$ and satisfying conditions of lemma \ref{jsympos}. 

Consider an arbitrary vertex $v_0\in F_I\subset P$. Denote by $r_1, \ldots, r_k$, $k=\dim F_I$, the edges of the face $F_I$ containing $v_0$.
If $\varphi_{a_{v_0,r_i}}\colon \Z_P\to \m C$ are the corresponding monomial maps, then there exists $\mu\in \m C^*$ such that $\varphi_{a_{v_0, r_i}}(z')$ = $\mu \varphi_{a_{v_0, r_i}}(z)$ for all $i\in [1,k]$. Consider the maps
$$
\phi_i(z) = \log \Bigl( \frac{|\varphi_{a_{v_0, r_i}}(z)|^2}{|\varphi_{b_{v_0}}(z)|^2} \Bigr), \quad i \in [1, k], \quad z \in \Z_P.
$$
Then $\phi_i(z')$ = $\phi_i(z)$ for all $i\in[1, k]$. 

Note that every function $\phi_i$ considered as the function of argument $w$ is linear. If ${\Phi = (\phi_1(w),~\ldots,~\phi_k(w))}$ is the corresponding linear map from $\m R^J$ to $\m R^k$, then $\Phi(w') = \Phi(w)$.
We need to show that  $\beta'$~=~$\beta$ or, equivalently, $w'$ = $w$.  
Since the map $j_F$ is injective, $w'$ = $w$ if and only if $x'$ and $x$ coincide.

It is enough to show that $\Phi(w')$ = $\Phi(w)$ implies that $x'$ = $x$ in~$j_F^{-1}(F_I)\subset \m R^k$. Denote by $p_J\colon \m R^m\to \m R^J$ the standard projection map. The composite map 
$\Phi \circ p_J \circ j_F\colon \m R^k\to \m R^k$ is a smooth (not linear) function of $x\in \m R^k$.
The map $\Phi \circ p_J$ is linear by $w$. The linear function $\phi_i\circ p_J$ is defined by the $m$-dimensional coefficient vector ${u_{r_i} = (a_{v_0, r_i} - b_{v_0})}$. By lemma \ref{samecharacter}, we have $u_{r_i} \in \op{Im}A_P$ = $\op{Im}\Lambda^T$. This means that the vectors $u_{r_i}$, which are also rows of the matrix $\Phi\circ p_J$, are linear combinations of columns of the matrix $A_P$.

Consider the matrix $L_1$ of a linear map corresponding to the affine embedding ${i_P\circ j_F\colon \m R^k\to \m R^m}$. The rank of $L_1$ is equal to $k$ and its rows are also the linear combinations of columns of $A_P$. Moreover, all columns of $L_1$ with indices from the set $I$ are zero.
Denote by $L_0$ the matrix of the mapping $\Phi \circ p_J$ which is linear by $w$ = $(w_j)$. As said before, the matrix $L_0$ is formed by vectors $u_{r_i}$. By corollary \ref{coordzero}, all coordinates of the vectors $u_{r_i}$, $i\in [1,k]$, from the index set $I$ are zero. By lemma \ref{unitmatrix}, the matrix $L_0$ contains the unit submatrix of rank $k$, so $L_0$ also has the rank $k$. The fact that the matrix $A_P^T = \Lambda$ is a characteristic function now implies that every row of $L_1$ is a linear combination of columns of $L_0$. This means that there exists a~$(k\times k)$-matrix $D$ such that $DL_0$ = $L_1^T$. Since $L_1$ has also rank $k$, the matrix $D$ is invertible. We can see that the matrix of composition of $\Phi \circ p_J$ with an automorphism of $\m R^k$ given by $D$ is exactly the matrix $L_1^T$. 

Consider now the composite map
$$
 U \stackrel{j_F}{\xrightarrow{\hspace*{0.5cm}}} \m R^m \stackrel{D\circ \Phi \circ p_J}{\xrightarrow{\hspace*{1.5cm}}} \m R^k,
$$
where $U$ = $j_F^{-1}(\stackrel{\circ}{F_I})$ is a convex region in $\m R^k$. If $j_F(x) = L_1(x) + d$, where $d = (d_j)$ and $L_1 = (l_{ji})$, then we have
$$
\beta_j = l_{ji} x_i + d_j,  
\quad
\phi_k = l_{jk} \log (\beta_j),
$$
and this implies
$$
\frac{\partial \phi_k}{\partial x_i} 
= 
\sum\limits_{j\in J} \frac{\partial \phi_k}{\partial \beta_j} \frac{\partial \beta_j}{\partial x_i}
=  
\sum\limits_{j\in J} \frac{l_{ji}}{\beta_j}  l_{jk}.
$$
We see that Jacobian of this composite map is symmetric and positive definite in every point of $U$. The proof of lemma \ref{modules} is accomplished by applying proposition \ref{jsympos}. The fact that Jacobian is non-degenerate in every point of $U$ also implies that $\tilde\varphi_{\m P, \tilde k_P}$ is a smooth embedding.
$\Box$

Note that the form $\omega = \sum\limits_{i=1}^n \phi_i dx_i$ in our case coincides with the form $i_P^* \Bigl(\sum\limits_{j\notin I} \log(\beta_j) d\beta_j\Bigr)$, which is obviously exact.

Let us consider an example of toric variety of complex dimension three over Stasheff polytope $K_5$.

The matrix $\Lambda = A_P^T$ has the form
$$
\Lambda = 
\begin{pmatrix}
1 & 0 & 0 & -1 & 0 & 0 & 0 & 1 & -1\\
0 & 1 & 0 & 0 & -1 & 0 & -1 & 0 & 1\\
0 & 0 & 1 & 0 & 0 & -1 & -1 & 1 & 0\\
\end{pmatrix},
$$
and the vector $b_P$ is equal to $(0, 0, 0, 3, 3, 3, 5, -1, 2)^T$. The matrix ${C = C_P^T}$ may be chosen in the form containing only zeros and ones:
$$
C = 
\begin{pmatrix}
1 & 0 & 0 & 1 & 0 & 0\\
0 & 1 & 0 & 0 & 0 & 1\\
0 & 0 & 1 & 1 & 0 & 1\\
1 & 0 & 0 & 0 & 0 & 0\\
0 & 1 & 0 & 0 & 0 & 0\\
0 & 0 & 1 & 0 & 0 & 0\\
0 & 0 & 0 & 1 & 1 & 1\\
0 & 0 & 0 & 0 & 1 & 0\\
0 & 0 & 0 & 1 & 1 & 0\\
\end{pmatrix}.
$$

Note that fixing $C = C_P^T$ immediately allows us to write explicitly all of the equations that define the moment-angle manifold ${\Z_P\subset \m C^9}$. When all of entries of $C$ are non-negative, the manifold $\Z_P$ is an~intersection of spheres and  cylinders in corresponding Euclidean space. In our case the manifold $\Z_P\subset\m C^9$ is given by the following equations:
$$
|z_1|^2+ |z_4|^2 = 3,
$$
$$
|z_2|^2 + |z_5|^2 = 3,
$$
$$
|z_3|^2 + |z_6|^2 = 3,
$$
$$
|z_1|^2 + |z_3|^2 + |z_7|^2 + |z_9|^2 = 7,
$$
$$
|z_7|^2 + |z_8|^2 + |z_9|^2 = 6,
$$
$$
|z_2|^2 + |z_3|^2 + |z_7|^2 = 5.
$$

The moment-angle manifold $\Z_P\subset \m C^9$ has a real dimension $12$. It is endowed with the action of the compact torus $\m T^9$, with six-dimensional subgroup $K\subset T^9$ acting freely on $\Z_P$. 

Stasheff polytope may be obtained from the cube $I^3$ by a sequence of cutting faces of codimension two. The polytopes that are obtained from the $n$-dimensional cube in this way are called {\it $2$-truncated cubes}. The theory of $2$-truncated cubes has been developed in the paper \cite{buchvol}; the class of $2$-truncated cubes provides a wide variety of quasitoric manifolds with an explicit description of the data $(P, \Lambda)$. When a face of codimension two of the simple polytope $P$ is cut, the corresponding matrix $\Lambda = A_P^T$ is extended by a column that is a sum of two other columns of $\Lambda$; these columns correspond to facets adjacent to the face being cut.

If we enumerate the facets of Stasheff polytope according to the rows of the matrix $C$ (see above), then two-element sets corresponding to the edges of $K_5$ are (1, 2), (1, 5), (1, 6), (1, 7), (1, 8), (2, 3), (2,~6), (2, 8), (2, 9), (3, 4), (3, 5), (3, 8), (3, 9), (4, 5), (4, 6), (4, 7), (4, 9), (5, 7), (5, 8), (6, 7), (6, 9). For every pair of numbers we need to exclude two corresponding rows from the matrix $C$ and then find the cokernel matrix (it would be a vector in $\m Z^7$). Let us write out all of the corresponding maps $\varphi_{ij}\colon \m C^9\to \m C$:\\

\begin{tabbing}
$\varphi_{12}~=~z_3\bar z_6 \bar z_7 z_8,$ \qquad\qquad\=
$\varphi_{15}~=~\bar z_3 z_6 z_7 \bar z_8,$ \qquad\qquad\=
$\varphi_{16}~=~\bar z_2 z_5 z_7 \bar z_9,$\\

$\varphi_{17}~=~\bar z_2 z_3 z_5 \bar z_6 z_8 \bar z_9,$\>
$\varphi_{18}~=~z_2 \bar z_5 \bar z_7 z_9,$\>
$\varphi_{23}~=~z_1 \bar z_4 z_8 \bar z_9,$\\

$
\varphi_{26}~=~z_1 \bar z_4 z_8 \bar z_9,
$\>
$
\varphi_{28}~=~\bar z_1 z_3 z_4 \bar z_6 \bar z_7 z_9,
$\>
$
\varphi_{29}~=~\bar z_3 z_6 z_7 \bar z_8,
$\\
$
\varphi_{34}~=~\bar z_2 z_5 z_7 \bar z_9,
$\>
$
\varphi_{35}~=~\bar z_1 z_4 \bar z_8 z_9,
$\>
$
\varphi_{38}~=~z_2 \bar z_5 \bar z_7 z_9,
$\\
$
\varphi_{39}~=~z_1 z_2 \bar z_4 \bar z_5 \bar z_7 z_8,
$\>
$
\varphi_{45}~=~z_3 \bar z_6 \bar z_7 z_8,
$\>
$
\varphi_{46}~=~z_2 \bar z_5 \bar z_7 z_9,
$\\
$
\varphi_{47}~=~\bar z_2 z_3 z_5 \bar z_6 z_8 \bar z_9,
$\>
$
\varphi_{57}~=~\bar z_1 z_4 \bar z_8 z_9,
$\>
$
\varphi_{58}~=~\bar z_1 z_3 z_4 \bar z_6 \bar z_7 z_9,
$\\
$
\varphi_{67}~=~\bar z_1 z_4 \bar z_8 z_9,
$\>
$
\varphi_{69}~=~z_1 z_2 \bar z_4 \bar z_5 \bar z_7 z_8.
$\\
\end{tabbing}

This example shows that even in non-trivial cases, when $P$ is not a~product of other polytopes, the number  $q$ of pairwise distinct vectors $u_r$ and pairwise distinct monomials $\varphi_{u_r}$ may be significantly less than $f_1(P)$ (the upper bound). In the case of Stasheff polytope we have $q = 6$ and $f_1(P)$ = $21$. These monomial functions embed the quasitoric manifold corresponding to Stasheff polytope to Euclidean space ${\m R^3\times \m C^6}$. 

Let us now construct the projective embedding of Stasheff polytope described in theorem \ref{projtoric}. According to the algorithm, we need to choose the set of lattice points in $\m Z^m$ formed by vertices of the polytope $i_P(P)\subset\m R^m_{\geqslant 0}$ and the lattice points on the edges that are next to vertices. These lattice points determine monomial functions which in turn define an embedding of the manifold.

The Stasheff polytope has $14$ vertices -- these are the points ${(1, 0, 0)}$, ${(0, 0, 1)}$, ${(1, 3, 0)}$, ${(0, 3, 1)}$, ${(2, 0, 0)}$, ${(3, 1, 0)}$, ${(3, 1, 3)}$, ${(2, 0, 3)}$, ${(3, 2, 3)}$, ${(0, 2, 3)}$, ${(0, 3, 2)}$, ${(3, 3, 2)}$, ${(0, 0, 3)}$, ${(3, 3, 0)}$. The additional points on the edges of $P$ are ${(1, 1, 0)}$, ${(1, 2, 0)}$, ${(2, 3, 0)}$, ${(3, 2, 0)}$, ${(0, 1, 1)}$, ${(0, 2, 1)}$, ${(3, 3, 1)}$, ${(3, 1, 1)}$, ${(2, 0, 1)}$, ${(1, 3, 2)}$, ${(2, 3, 2)}$, ${(3, 1, 2)}$, ${(2, 0, 2)}$, ${(0, 0, 2)}$, ${(1, 0, 3)}$, ${(0, 1, 3)}$, ${(1, 2, 3)}$, ${(2, 2, 3)}$ -- $18$ points in total.

The following monomial functions define a projective embedding of the manifold corresponding to Stasheff polytope $K_5$:\\

\begin{tabbing}
$\varphi_{b_1} = z_1 z_4^2 z_5^3 z_6^3 z_7^5 z_9,$ \qquad\qquad\=
$\varphi_{b_2} = z_3 z_4^3 z_5^3 z_6^2 z_7^4 z_9^2,$ \qquad\qquad\=
$\varphi_{b_3} = z_1 z_2^3 z_4^2 z_6^3 z_7^2 z_9^4,$\\
$\varphi_{b_4} = z_2^3 z_3 z_4^3 z_6^2 z_7 z_9^5,$\>
$\varphi_{b_5} = z_1^2 z_4 z_5^3 z_6^3 z_7^5 z_8,$\>
$\varphi_{b_6} = z_1^3 z_2 z_5^2 z_6^3 z_7^4 z_8^2,$\\
$\varphi_{b_7} = z_1^3 z_2 z_3^3 z_5^2 z_7 z_8^5,$\>
$\varphi_{b_8} = z_1^2 z_3^3 z_4 z_5^3 z_7^2 z_8^4,$\>
$\varphi_{b_9} = z_1^3 z_2^2 z_3^3 z_5 z_8^5 z_9,$\\  
$\varphi_{b_{10}} = z_2^2 z_3^3 z_4^3 z_5 z_8^2 z_9^4,$\>
$\varphi_{b_{11}} = z_2^3 z_3^2 z_4^3 z_6 z_8 z_9^5,$\>
$\varphi_{b_{12}} = z_1^3 z_2^3 z_3^2 z_6 z_8^4 z_9^2,$\\
$\varphi_{b_{13}} = z_3^3 z_4^3 z_5^3 z_7^2 z_8^2 z_9^2,$\>
$\varphi_{b_{14}} = z_1^3 z_2^3 z_6^3 z_7^2 z_8^2 z_9^2.$
\end{tabbing}

\vspace{0.5cm}

\begin{tabbing}
$\varphi_{a_1} = z_1 z_2 z_4^2 z_5^2 z_6^3 z_7^4 z_9^2,$
\qquad\qquad\=
$\varphi_{a_2} = z_1 z_2^2 z_4^2 z_5 z_6^3 z_7^3 z_9^3,$
\qquad\qquad\=
$\varphi_{a_3} = z_1^2 z_2^3 z_4 z_6^3 z_7^2 z_8 z_9^3,$\\
$\varphi_{a_4} = z_1^3 z_2^2 z_5 z_6^3 z_7^3 z_8^2 z_9,$\>
$\varphi_{a_5} = z_2 z_3 z_4^3 z_5^2 z_6^2 z_7^3 z_9^3,$\>
$\varphi_{a_6} = z_2^2 z_3 z_4^3 z_5 z_6^2 z_7^2 z_9^4,$\\
$\varphi_{a_7} = z_1^3 z_2^3 z_3^1 z_6^2 z_7 z_8^3 z_9^2,$\>
$\varphi_{a_8} = z_1^3 z_2 z_3 z_5^2 z_6^2 z_7^3 z_8^3,$\>
$\varphi_{a_9} = z_1^2 z_3 z_4 z_5^3 z_6^2 z_7^4 z_8^2,$\\
$\varphi_{a_{10}} = z_1 z_2^3 z_3^2 z_4^2 z_6 z_8^2 z_9^4,$\>
$\varphi_{a_{11}} = z_1^2 z_2^3 z_3^2 z_4 z_6 z_8^3 z_9^3,$\>
$\varphi_{a_{12}} = z_1^3 z_2 z_3^2 z_5^2 z_6 z_7^2 z_8^4,$\\
$\varphi_{a_{13}} = z_1^2 z_3^2 z_4 z_5^3 z_6 z_7^3 z_8^3,$\>
$\varphi_{a_{14}} = z_3^2 z_4^3 z_5^3 z_6 z_7^3 z_8 z_9^2,$\>
$\varphi_{a_{15}} = z_1 z_3^3 z_4^2 z_5^3 z_7^2 z_8^3 z_9,$\\
$\varphi_{a_{16}} = z_2 z_3^3 z_4^3 z_5^2 z_7 z_8^2 z_9^3,$\>
$\varphi_{a_{17}} = z_1 z_2^2 z_3^3 z_4^2 z_5 z_8^3 z_9^3,$\>
$\varphi_{a_{18}} = z_1^2 z_2^2 z_3^3 z_4 z_5 z_8^4 z_9^2.$\\
\end{tabbing}

We see that the dimension of the projective embedding is much more than the dimension of the affine embedding, so in this sense the case of Stasheff polytope $K_5$ is quite different from the manifold $\m CP^n$ corresponding to the simplex $\Delta^n$.

\end{document}